\documentclass[11pt]{article}
\usepackage{amssymb}
\usepackage{amsfonts}
\usepackage{amsmath}
\usepackage{amsthm}
\usepackage{graphicx}
\usepackage{epsfig}
\usepackage{psfrag}
\usepackage{color}
\usepackage{dsfont}
\makeatletter
\xdef\@endgadget#1{{\unskip\nobreak\hfil\penalty50\hskip1em\hbox{}\nobreak
    \hfil#1\parfillskip=0pt\finalhyphendemerits=0\par}}
\def\@qedsymbol{${}_\blacksquare$}
\def\qed{\@endgadget{\@qedsymbol}}
\newtheorem{lemma}{Lemma}[section]

\newtheorem{thm}[lemma]{Theorem}

\newtheorem{example}[lemma]{Example}

\newtheorem{defn}[lemma]{Definition}

\newtheorem{prop}[lemma]{Proposition}
\newtheorem{remark}[lemma]{Remark}

\def\Real{\mathbb{R}}

      \def\HH{\mathcal{H}}

\newcommand{\Ltwo}{\boldsymbol{\rm L}_{2}}
\newcommand{\Ltwoe}{\boldsymbol{\rm L}_{2e}}

\title{The converse of the passivity and small-gain theorems for input-output maps}

\author{Sei Zhen Khong, Arjan van der Schaft
\thanks{Sei Zhen Khong is with the Institute for Mathematics and its Applications, The University of Minnesota, Minneapolis, MN 55455, USA, {\tt\small szkhong@umn.edu}, 
Arjan van der Schaft is with the Johann Bernoulli Institute for Mathematics and Computer Science, University of Groningen, Groningen, the Netherlands, {\tt\small a.j.van.der.schaft@rug.nl}}
}
\date{Version: June 25, 2018; accepted for publication in {\it Automatica}}
\begin{document}

\maketitle
\thispagestyle{empty}
\pagestyle{empty}
\begin{abstract}
  We prove the following converse of the passivity theorem. Consider a causal system given by a sum of a linear time-invariant and a passive linear
  time-varying input-output map. Then, in order to guarantee stability (in the sense of finite $\Ltwo$-gain) of the feedback interconnection of the
  system with an {\it arbitrary} nonlinear output strictly passive system, the given system must itself be output strictly passive. The proof is based
  on the S-procedure lossless theorem. We discuss the importance of this result for the control of systems interacting with an output strictly
  passive, but otherwise completely unknown, environment. Similarly, we prove the necessity of the small-gain condition for closed-loop stability of certain
  time-varying systems, extending the well-known necessity result in linear robust control.
\end{abstract}


\section{Introduction}
The {\it passivity} and {\it small-gain} theorems are fundamental to large parts of systems and control theory, see
e.g. \cite{Wil72, MoyHil78, Vid81, MegRan97, Sch17}. Both theorems provide a stability `certificate' when feedback interconnecting the given system with an arbitrary system which is either (in the small-gain setting) assumed to have an $L_2$-gain smaller than the reciprocal of the $L_2$-gain of the given system, or is (output strictly) passive like the given system. These theorems are valid from linear finite-dimensional systems to nonlinear and infinite-dimensional systems.

The current paper is concerned with the {\it converse} of these theorems; that is the {\it necessity} of the (strict) passivity or the small-gain
condition for closed-loop stability when feedback interconnecting a given system with an {\it arbitrary} system, which is {\it unknown} apart from a
passivity or $\Ltwo$-gain assumption. Surprisingly, this converse of the {\it passivity} theorem has hardly been studied in the literature; despite its
fundamental importance in applications. For example, in order to guarantee stability of a controlled robotic system interacting with a passive, but else completely {\it unknown}, environment, the converse of the passivity theorem tells us that the controlled robot {\it must} be output strictly passive as seen from the interaction port of the robot with the environment. This has far-reaching methodological implications for control
design, since it means that rendering by control the system output strictly passive at the interaction port is not only a valid option, but is the {\it only} option guaranteeing stability for an unknown passive environment. The same holds within the context of robust nonlinear
control whenever we replace `environment' by the uncertain part of the system.

Up to now this converse passivity theorem was only proved for {\it linear time-invariant single-input single-output} systems in \cite{colgate-hogan},
using arguments from Nyquist stability theory\footnote{Roughly speaking, by showing that if $\Sigma_1$ is {\it not} passive, a positive-real transfer
  function (corresponding to a passive system $\Sigma_2$) can be constructed such that the closed-loop system fails the Nyquist stability test.},
exactly with the robotics motivation in mind. The same motivation was elaborated on in \cite{Str15}, where the following form of a converse passivity
theorem was obtained for nonlinear systems in state space form. If a system is {\it not} passive then for any given constant $K$ one can define a
passive system that extracts from the given system an amount of energy that is larger than $K$, implying that the norm of the state of the constructed
system becomes larger than $K$, thereby demonstrating some sort of instability of the closed-loop system.
In the present paper, a converse of the passivity theorem will be derived for a class of {\it input-output maps}, namely those decomposable
into a sum of a linear time-invariant map and a passive linear time-varying map. This converse passivity theorem involves feedback interconnections
with nonlinear systems and will be formulated in three versions in Section~\ref{sec: passivity}, with their own range of applicability. In all cases
the proofs are based on the S-procedure lossless theorem due to Megretski \& Treil \cite{MegTre93}; see also \cite[Thm. 7]{Jon01}.

Converse statements of the {\it small-gain} theorem are much more present in the literature; see e.g. \cite[Theorem 9.1]{ZDG} for the
finite-dimensional linear case and \cite{CZ} for infinite-dimensional linear systems. However, to the best of our knowledge, the converse of the small-gain theorem for {\it linear time-varying} systems interconnected in feedback with nonlinear systems, as obtained in Section 4, is new, while also the proof line is different from the existing one. Similarly to the passivity case, this converse will be formulated for a class of linear time-varying input-output maps, and the proofs, in two
different
versions, will be based on the S-procedure lossless theorem. \\
Finally, Section 5 presents the conclusions, and discusses problems for further research.
A preliminary version of some of the results in Section 3 of this paper was presented at the IFAC World Congress 2017; cf. \cite{Kho-Sch16}.

\section{Preliminaries}

This section summarizes the background for this paper; see e.g. \cite{Sch17} for details.
Denote the set of $\Real^n$-valued Lebesgue square-integrable functions by
\[
\Ltwo^n := \Big\{v : [0, \infty) \to \Real^n \mid  \|v\|_2^2 := \int_0^\infty v(t)^Tv(t) \, dt < \infty\Big\}.
\]
For any two $v,w \in \Ltwo^n$ denote the $\Ltwo^n$-inner product as
\[
\langle v, w \rangle := \int_0^\infty v(t)^Tw(t) \, dt
\]
Define the {\it truncation} operator $(P_Tv)(t) := v(t)$ for $t \leq T$; $(P_Tv)(t) := 0$ for $t > T$, and the extended function space
\begin{align*}
 \Ltwoe^n & := \{v : [0, \infty) \to \Real^n \mid P_Tv \in \Ltwo, \; \forall T \in [0, \infty)\}.
\end{align*}
In what follows, the superscript $n$ will often be suppressed for notational simplicity.
Throughout this paper a {\it system} will be specified by an input-output map $\Delta : \Ltwoe^m \to \Ltwoe^p$ satisfying $\Delta(0)=0$. \\
Define for any $\tau \geq 0$ the {\it right shift} operator $(S_{\tau} (u))(t) = u(t - \tau)$ for $t\geq \tau$ and $(S_{\tau} (u))(t) =0$ for $0\leq t <\tau$. The system $\Delta$ is said to be \emph{time-invariant} if $S_{\tau} \Delta = \Delta S_{\tau}$ for all $\tau > 0$. 
Furthermore, the system $\Delta$ is {\it bounded} if $\Delta$ maps $\Ltwo^m$ into $\Ltwo^p$. It is said to have {\it $\Ltwo$-gain} $\leq \gamma$ for some $\gamma > 0$ ({\it finite} $\Ltwo$-gain) if
\begin{equation}\label{eq: finitegain}
\| P_T \Delta (u) \|_2 \leq \gamma  \| P_T u\|_2
\end{equation}
for all $u \in \Ltwoe^m$ and all $T\geq 0$. The infimum of all $\gamma$ satisfying \eqref{eq: finitegain} is called the $\Ltwo$-gain of $\Delta$. 
The system $\Delta$ is \emph{causal} if $P_T \Delta P_T = P_T\Delta$ for all $T \geq 0$. It is well-known, see e.g. \cite[Proposition 1.2.3]{Sch17}, that a causal system $\Delta$ has finite $\Ltwo$-gain if and only if, instead of \eqref{eq: finitegain},
\begin{equation}\label{eq: finitegain1}
\|\Delta (u) \|_2 \leq \gamma  \| u\|_2
\end{equation}
for all $u \in \Ltwo^m$.  For the purpose of {\it interconnection} of systems the above notions are generalized from maps to {\it relations}
$R \subset \Ltwoe^m \times \Ltwoe^p$ satisfying $(0,0) \in R$ as follows \cite{Sch17}. A relation $R$ is said to be {\it bounded} if whenever
$(u,y)\in R$ and $u\in \Ltwo$ then also $y\in \Ltwo$. Furthermore, $R$ has {\it finite $\Ltwo$-gain} if
\begin{equation}\label{eq: finitegain2}
\| P_Ty\|_2 \leq \gamma  \| P_T u\|_2
\end{equation}
for all $(u,y) \in R$ and all $T\geq 0$. Also, $R$ is said to be {\it causal} if whenever $(u_1,y_1) \in R$, $(u_2,y_2)\in R$ satisfy $P_Tu_1=P_Tu_2$, then $P_Ty_1=P_Ty_2$. A causal relation $R$ has finite $\Ltwo$-gain if and only i,f instead of \eqref{eq: finitegain2} ,
\begin{equation}\label{eq: finitegain3}
\|y\|_2 \leq \gamma  \| u\|_2
\end{equation}
for all $\Ltwo$ pairs $(u,y) \in R$.
The system $\Delta:   \Ltwoe^m \to \Ltwoe^m$ (i.e., $p=m$) is said to be \emph{passive}~\cite{Wil72, Vid81} if
\begin{align} \label{eq: passive}
\int_0^T u(t)^T (\Delta (u))(t) \, dt \geq 0, 
\end{align}
for all $u \in \Ltwoe, T >0$. Furthermore, it is called \emph{strictly passive} if there exist $\delta>0, \epsilon > 0$ such that
\[
\int_0^T u(t)^T (\Delta (u))(t) \, dt \geq \delta \|P_T u\|_2^2 +  \epsilon \|P_T \Delta (u)\|_2^2
\]
for all $u \in \Ltwoe, T >0$, and \emph{output strictly passive} if this holds with $\delta=0$.
In case $\Delta$ is bounded and causal, then passivity is equivalent~\cite[Proposition 2.2.5]{Sch17} to
\begin{align} \label{eq: causal_passive}
\int_0^\infty u(t)^T (\Delta (u))(t) \, dt \geq 0
\end{align}
for all $u \in \Ltwo^m$. (Note that the integral is well-defined because of boundedness of $\Delta$ and the Cauchy-Schwartz inequality.)
Similarly, in this case $\Delta$ is {\it strictly passive} if there exist $\delta>0, \epsilon > 0$ such that
\begin{equation}\label{eq: causal_passive1}
\int_0^\infty u(t)^T (\Delta (u))(t) \, dt \geq \delta \|u\|_2^2 + \epsilon \|\Delta (u)\|_2^2 \quad \forall u \in \Ltwo^m,
\end{equation}
and {\it output strictly passive} if this holds with $\delta=0$.
For later use we also recall the basic property that any output strictly passive system has finite $\Ltwo$-gain; cf. \cite[Theorem 2.2.13]{Sch17}.
Like in the $\Ltwo$-case these passivity notions are directly extended to {\it relations} $R \subset \Ltwoe^m \times \Ltwoe^m$ satisfying
$(0,0) \in R$. Indeed, $R$ is called {\it strictly passive} if there exist $\delta>0, \epsilon > 0$ such that for all\footnote{Throughout it us assumed that all integrals are well-defined.} $(u,y) \in R$, $T > 0$
\begin{equation}
\int_0^T u(t)^T y(t) \, dt \geq \delta \|P_T u\|_2^2 +  \epsilon \|P_T y\|_2^2,
\end{equation}
and {\it output strictly passive} if this holds with $\delta=0$. Furthermore, a bounded causal relation $R$ is strictly passive if there exist
$\delta>0, \epsilon > 0$ such that for all $(u,y) \in R$
\begin{equation}
\int_0^{\infty} u(t)^T y(t) \, dt \geq \delta \|u\|_2^2 +  \epsilon \|y\|_2^2,
\end{equation}
and output strictly passive if this holds with $\delta=0$. 
%

\setlength{\unitlength}{0.7cm}
\begin{figure}[h]
  \centering 
  \includegraphics[scale=0.6]{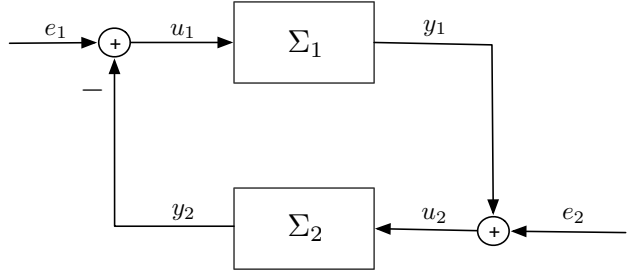}
  \caption{Feedback configuration} \label{fig: feedback}
\end{figure}

The main object of study in this paper is the {\it feedback interconnection} of two systems $\Sigma_1 : \Ltwoe^{m_1} \to \Ltwoe^{p_1}$ and $\Sigma_2 : \Ltwoe^{m_2} \to \Ltwoe^{p_2}$, with $m_1=p_2, m_2=p_1$, described by (see Figure~\ref{fig: feedback})
\begin{equation} \label{eq: FB}
\begin{array}{rcllll}
 u_1 & = & e_1 - y_2, \quad & u_2 &= &e_2 + y_1, \\[2mm]
 y_1 & = & \Sigma_1(u_1), \quad & y_2 &= &\Sigma_2(u_2).
\end{array}
\end{equation}
The resulting {\it closed-loop system} with inputs $(e_1,e_2)$ and outputs $(y_1,y_2)$ will be denoted by $\Sigma_1 \| \Sigma_2$, and defines by \eqref{eq: FB} a relation in the space of all $(e_1,e_2,u_1,u_2,y_1,y_2) \in \Ltwoe$. Projection on the space of $(e_1,e_2,u_1,u_2) \in \Ltwoe$, respectively of $(e_1,e_2,y_1,y_2) \in \Ltwoe$, yields the relations
\[
\begin{array}{rcl}
R_{eu} \! &:= & \! \{(e_1,e_2,u_1,u_2) \in \Ltwoe \! \mid \! \eqref{eq: FB} \mbox{ holds for some } y_1,y_2\} \\[2mm]
R_{ey} \! &:= & \! \{(e_1,e_2,y_1,y_2) \in \Ltwoe \! \mid \! \eqref{eq: FB} \mbox{ holds for some } u_1,u_2\}
\end{array}
\]

\begin{defn} \label{def: FB} The feedback interconnection $\Sigma_1 \| \Sigma_2$ has \emph{finite $\Ltwo$-gain} if the relation $R_{eu}$, or equivalently (see \cite[Lemma 1.2.12]{Sch17}, the relation $R_{ey}$, has finite $\Ltwo$-gain. 
$\Sigma_1 \| \Sigma_2$ with $m_1=m_2=p_1=p_2$ is said to be \emph{passive} whenever the relation $R_{ey}$ is passive. 
The feedback interconnection for $e_2 = 0$, denoted by $\Sigma_1 \|_{e_2=0} \Sigma_2$, is said to have finite $\Ltwo$-gain if the corresponding relation $R_{e_1y_1}$ has finite $\Ltwo$-gain, 
and is said to be passive if $R_{e_1y_1}$ is passive. For notational convenience, we denote the map from $e_1$ to $y_1$ by $\Sigma_1 \|_{e_2=0} \Sigma_2$.
\end{defn}
Finally, if the systems $\Sigma_1$ and $\Sigma_2$ are causal, then so are the relations $R_{ey}$ and $R_{eu}$; see \cite[Proposition
1.2.14]{Sch17}. The same statement is easily seen to hold for $R_{e_1y_1}$. All systems and relations are taken to be causal throughout this paper.
%


%

\section{Passivity as a necessary condition for stable interaction} \label{sec: passivity}

The classical passivity theorem, see e.g. \cite{Sch17}, asserts that the feedback interconnection $\Sigma_1 \| \Sigma_2$ of two passive systems $\Sigma_1, \Sigma_2$ is again a passive system. Similarly, the interconnected system $\Sigma_1 \|_{e_2=0} \Sigma_2$ is passive. 
In this section we will derive a converse passivity theorem\footnote{A {\it different}, and easy to prove, converse result stating that passivity of $\Sigma_1 \| \Sigma_2$ implies that both $\Sigma_1$ and $\Sigma_2$ are passive was formulated in \cite{kerber-vds11}, \cite{Kho-Sch16}; see also \cite[Proposition 4.3.8]{Sch17}.} stating that a necessary condition in order that
any closed-loop system arising from interconnecting a {\it given} system $\Sigma_1$ to with an {\it unknown, but output strictly passive}, system $\Sigma_2$ is {\it stable} (in the sense of having finite $\Ltwo$-gain), is that the system $\Sigma_1$ is itself output strictly passive. 
As already indicated in the introduction, this result is crucial e.g. in the control of robotic systems; see also the discussion and example at the end of this section.  We will formulate
three different versions of this theorem. Before doing so we first state the following version of the S-procedure lossless
theorem, which can be obtained from~\cite[Thm. 7 and Ex. 28]{Jon01}, based on \cite{MegTre93}.

\begin{prop}[S-procedure lossless theorem] \label{prop: Sprod} Let $\HH \subset \Ltwo$ be a vector space satisfying $S_\tau \HH \subset \HH$ for all
  $\tau \geq 0$ and $\sigma_i : \HH \to \Real$ be defined as $\sigma_i(f) := \langle f, \Phi_i f \rangle$, where $\Phi_i = \Phi_i^T$ is a constant
  matrix, $i = 0, 1$. Suppose there exists an $f^* \in \HH$ such that $\sigma_1(f^*) > 0$, then the following are equivalent:
\begin{enumerate} \renewcommand{\theenumi}{\textup{(\roman{enumi})}}\renewcommand{\labelenumi}{\theenumi}
 \item $\sigma_0(f) \leq 0$ for all $f \in \HH$ such that $\sigma_1(f) \geq 0$;

 \item $\exists \mu \geq 0$ such that $\sigma_0(f) + \mu \sigma_1(f) \leq 0, \quad \forall f \in \HH$.
\end{enumerate}
\end{prop}

The first version of the converse passivity theorem is as follows.

\begin{thm} \label{thm: passive} Given bounded $\Sigma_1 = G + \Delta$, where $G$ is linear time-invariant and $\Delta$ is linear passive, then there
  exists $\gamma > 0$ such that the closed-loop system $\Sigma_1 \| \Sigma_2$ has $\Ltwo$-gain $\leq \gamma$ for {\it all} bounded passive $\Sigma_2$
  if and only if $\Sigma_1$ is strictly passive.
\end{thm}

\noindent
{\bf Proof} \,
{\it Sufficiency} is well known in the literature. Indeed, strict passivity of $\Sigma_1$ together with passivity of $\Sigma_2$ yields
\begin{align*}
\epsilon (\|y_1\|_2^2 + \|u_1\|_2^2) & \leq \langle u_1, y_1 \rangle + \langle u_2, y_2 \rangle \\
& = \langle e_1 - y_2, y_1 \rangle + \langle e_2 + y_1, y_2 \rangle \\
& = \langle e_1, y_1 \rangle + \langle e_2, y_2 \rangle.
\end{align*}
Therefore, substituting $u_1=e_1-y_2, u_2=e_2+y_1$,
\[
\|y_1\|_2^2 + \langle e_1 - y_2, e_1 - y_2 \rangle \leq \frac{1}{\epsilon} (\langle e_1, y_1 \rangle + \langle e_2, y_2 \rangle),
\]
or
\[
\|y_1\|_2^2 + \|y_2\|_2^2 -2 \langle e_1, y_2 \rangle + \|e_1\|_2^2 \leq \frac{1}{\epsilon} (\langle e_1, y_1 \rangle + \langle e_2, y_2 \rangle).
\]
It follows that
\[
\|y\|_2^2 \leq 2 \langle e_1, y_2 \rangle + \frac{1}{\epsilon} \langle e, y \rangle \leq \left(2 + \frac{1}{\epsilon}\right) \|e\|_2\|y\|_2,
\]
where $y := (y_1, y_2)^T$ and $e := (e_1, e_2)^T$, and the Cauchy-Schwarz inequality has been used. Dividing both sides by $\|y\|_2$ the result follows.\\
To show {\it necessity}, let $\tilde{G}_\epsilon := G - \epsilon I$ and $\tilde{\Delta}_\epsilon := \Delta + \epsilon I$ for $\epsilon > 0$ so that
$\Sigma_1 = \tilde{G}_\epsilon + \tilde{\Delta}_\epsilon$. Recall from the theory of loop transformations~\cite[Section 3.5]{GreLim95} that the
finite $\Ltwo$-gain of $\Sigma_1 \| \Sigma_2$ is equivalent to that of $\tilde{G}_\epsilon\| [\Sigma_2 \|_{e_2 = 0} \tilde{\Delta}_\epsilon]$.
Furthermore, since $\tilde{\Delta}_\epsilon$ is strictly passive, that $\tilde{G}_\epsilon \| [\Sigma_2 \|_{e_2 = 0} \tilde{\Delta}_\epsilon]$ has
$\Ltwo$-gain $\leq \gamma$ for all bounded passive $\Sigma_2$ and $\epsilon > 0$ is equivalent to $G \| \Sigma_2$ having $\Ltwo$-gain $\leq \gamma$
for all bounded passive $\Sigma_2$. Define the vector space
\[
\HH := \{(u_1, u_2, e_1, e_2) \in \Ltwo\ |\ u_2 = e_2 + G (u_1)\}.
\]
Note that $S_\tau \HH \subset \HH$ for all $\tau \geq 0$ due to the time-invariance of $G$. 
Define now the following quadratic forms $\sigma_i : \HH
\to \Real$, $i = 0, 1$, as 
\begin{align*}
 \sigma_0(u_1, u_2, e_1, e_2) & := \left\langle 
\begin{bmatrix} 
 u_1 \\
 u_2 \\
 e_1 \\
 e_2
\end{bmatrix},
\begin{bmatrix}
I & 0 & 0 & 0 \\
0 & I & 0 & 0 \\
0 & 0 & -\gamma^2 I & 0 \\
0 & 0 & 0 & -\gamma^2 I
\end{bmatrix}
\begin{bmatrix} 
 u_1 \\
 u_2 \\
 e_1 \\
 e_2
\end{bmatrix}
\right\rangle 
\end{align*}
\begin{align*}
 \sigma_1(u_1, u_2, e_1, e_2) & := \frac{1}{2} \left\langle 
\begin{bmatrix} 
 u_1 \\
 u_2 \\
 e_1 \\
 e_2
\end{bmatrix}, 
\begin{bmatrix}
0 & -I & 0 & 0 \\
-I & 0 & I & 0 \\
0 & I &0 & 0 \\
0 & 0 & 0 & 0
\end{bmatrix}
\begin{bmatrix} 
 u_1 \\
 u_2 \\
 e_1 \\
 e_2
\end{bmatrix}
\right\rangle.
\end{align*}
Note that $\sigma_1(u_1, u_2, e_1, e_2) = u_2^T(e_1 - u_1)$, and hence it is easy to see that there exists $(u_1^*, u_2^*, e_1^*, e_2^*) \in \HH$ such
that $\sigma_1(u_1^*, u_2^*, e_1^*, e_2^*) > 0$.
It is immediately seen that $\sigma_0 \leq 0$ corresponds to the $\Ltwo$-gain of $R_{eu}$ being $ \leq \gamma $, while $\sigma_1 \geq 0$ corresponds
to any bounded  passive $\Sigma_2$. Hence, if the closed-loop system $G \| \Sigma_2$ has $\Ltwo$-gain $\leq \gamma$ for all
bounded  passive $\Sigma_2$, then
\begin{align*}
 \sigma_0(u_1, u_2, e_1, e_2) \leq 0 \quad \forall & (u_1, u_2, e_1, e_2) \in \HH \\
 & \text{such that} \; \sigma_1(u_1, u_2, e_1, e_2) \geq 0.
\end{align*}
This is equivalent, via the S-procedure lossless theorem (cf. Proposition~\ref{prop: Sprod}), to the existence of $\mu \geq 0$
such that
\begin{align*}
 \sigma_0(u_1, u_2, e_1, e_2) + \mu \sigma_1(u_1, u_2, & e_1, e_2) \leq 0, \\
 & \forall (u_1, u_2, e_1, e_2) \in \HH.
\end{align*}
Within the subset $\{(u_1, u_2, 0, 0) \in \Ltwo\ |\ u_2 = G (u_1)\} \subset \HH$, this yields
\[
\|G (u_1)\|_2^2 + \|u_1\|_2^2 - \mu \langle u_1, G (u_1) \rangle \leq 0, \quad \forall u_1 \in \Ltwo.
\]
This implies $\mu>0$, and thus
\[
\langle u_1, G (u_1) \rangle \geq \frac{1}{\mu} (\|G (u_1)\|_2^2 + \|u_1\|_2^2), \quad \forall u_1 \in \Ltwo,
\]
i.e., $G$ is strictly passive. Consequently, $\Sigma_1$ is strictly passive.

Roughly speaking, the new `only if' direction of the above theorem can be summarized by saying that in order that $\Sigma_1 \| \Sigma_2$ is {\it stable} (in the sense of having uniformly bounded $\Ltwo$-gain) for all passive $\Sigma_2$, then $\Sigma_1$ needs to be strictly passive. 
On the other hand, often in physical system examples (e.g., most mechanical systems) {\it output} strict passivity is a more natural property, since strict passivity can only occur for systems with direct feedthrough term;  see \cite[Proposition 4.1.2]{Sch17}. The following second version of the converse passivity theorem obviates this problem.

\setlength{\unitlength}{0.7cm}
\begin{figure}[h]
  \centering 
  \includegraphics[scale=0.5]{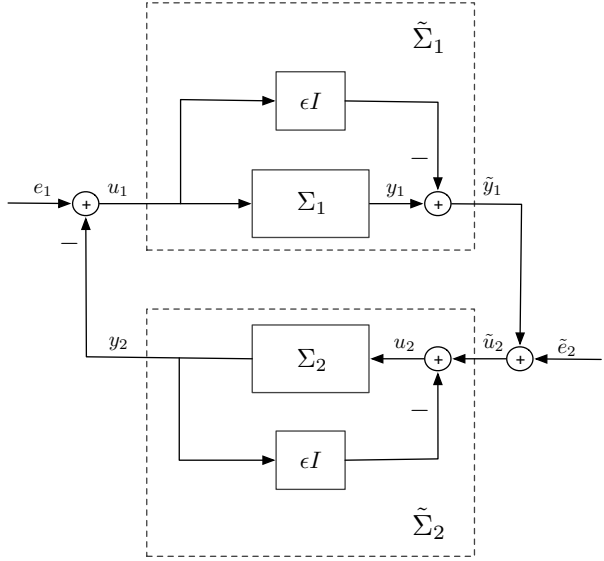}
  \caption{Loop transformation} \label{fig: loop_trans}
\end{figure}

\begin{thm}
\label{thm: passive2}
Given bounded  $\widetilde{\Sigma}_1 = G + \Delta$, where $G$ is linear time-invariant and $\Delta$ is linear passive, then
there exists $\gamma >0$ such that the closed-loop system $\widetilde{\Sigma}_1 \| \widetilde{\Sigma}_2$ has $\Ltwo$-gain $\leq \gamma$ for {\it all}
 output strictly passive $\widetilde{\Sigma}_2$ if and only if $\widetilde{\Sigma}_1$ is output strictly passive.
\end{thm}

\noindent
{\bf Proof} \,
  Sufficiency can be shown in a similar manner using the arguments in the sufficiency proof for Theorem~\ref{thm:
    passive}; see also \cite{Sch17}. For necessity, note that for any $\epsilon > 0$ the output strictly passive $\widetilde{\Sigma}_2$ can be written as the feedback interconnection of a bounded passive $\Sigma_2$ and $\epsilon I$, where $I$ denotes the identity operator. To see this, define $\Sigma_2$ as in Figure~\ref{fig: loop_trans}. Then by output strict passivity of $\widetilde{\Sigma}_2$
\[
\int_0^\infty \tilde{u}_2(t)^T y_2(t) \, dt \geq \epsilon \|y_2\|_2^2 \quad \forall \tilde{u}_2 \in \Ltwo,
\]
implying that
\[
\int_0^\infty (\tilde{u}_2(t) - \epsilon y_2(t))^T y_2(t) \, dt = \int_0^\infty u_2(t)^T y_2(t) \, dt \geq 0.
\]
The last inequality holds for all $u_2 \in \Ltwo$, since given any $u_2 \in \Ltwo$, $\tilde{u}_2 := (I + \epsilon \Sigma_2)u_2 \in \Ltwo$ yields the
desired $u_2$. It follows that $\Sigma_2$ is bounded and passive. By the same token, the negative feedback interconnection of a bounded passive
$\Sigma_2$ and $\epsilon I$ with $\epsilon > 0$ is output strictly passive. 
By defining $\Sigma_1 := \widetilde{\Sigma}_1 + \epsilon I$ as illustrated in Figure~\ref{fig: loop_trans}, one obtains the loop transformation
configuration therein. Since finite $\Ltwo$-gain of the closed-loop system $\widetilde{\Sigma}_1 \| \widetilde{\Sigma}_2$ in
Figure~\ref{fig: loop_trans} is equivalent to that of $\Sigma_1 \| \Sigma_2$ in Figure~\ref{fig: feedback}~\cite[Section 3.5]{GreLim95}, application
of Theorem~\ref{thm: passive} then yields that $\Sigma_1$ is strictly passive. For sufficiently small $\epsilon > 0$, it follows that
$\widetilde{\Sigma}_1 = \Sigma_1 - \epsilon I$ is output strictly passive. 

\medskip

Both Theorems \ref{thm: passive} and \ref{thm: passive2} require an exogenous signal $e_2$, which is often not the typical case in applications. This
motivates the following {\it third} version of the converse passivity theorem.

\begin{thm}\label{Theorem6}
  Given bounded $\Sigma_1 = G + \Delta$, where $G$ is linear time-invariant and $\Delta$ is linear passive, then there exists $\gamma >0$ such that
  the closed-loop system $\Sigma_1 \|_{e_2=0} \Sigma_2$ has $\Ltwo$-gain from $e_1$ to $y_1$ less than or equal to $\gamma$ for {\it all} bounded
  passive $\Sigma_2$ if and only if $\Sigma_1$ is output strictly passive.
\end{thm}

\noindent
{\bf Proof} \,
Sufficiency is well known in the literature. Indeed, if $\Sigma_1$ is output strictly passive and $\Sigma_2$ is passive, then for some $\varepsilon > 0$
 \begin{align*}
 \langle e_1, y_1 \rangle & = \langle u_1 + y_2, y_1 \rangle = \langle u_1, y_1 \rangle + \langle y_2, y_1 \rangle \\
 & = \langle u_1, y_1 \rangle + \langle u_2, y_2 \rangle \geq \varepsilon \| y_1\|_2^2,
 \end{align*}
 showing that the closed-loop system is $\varepsilon$-output strictly passive, and hence (see e.g. \cite{Sch17}, Theorem 2.2.13) has $\Ltwo$-gain $\leq \frac{1}{\varepsilon}$. 
To show necessity, note that by the same arguments in Theorem~\ref{thm: passive}, the hypothesis is equivalent
to $G \|_{e_2=0} \Sigma_2$ having $\Ltwo$-gain $\leq \gamma$ for all bounded  passive $\Sigma_2$. Define
\[
\HH := \{(u_1, y_1, e_1) \in \Ltwo\ |\ y_1 = G (u_1)\}
\]
and the quadratic forms $\sigma_i : \HH
\to \Real$, $i = 0, 1$, as 
\begin{align*}
 \sigma_0(u_1, y_1, e_1) & := \left\langle 
\begin{bmatrix} 
 u_1 \\
 y_1 \\
 e_1 
\end{bmatrix},
\begin{bmatrix}
0 & 0 & 0  \\
0 & I & 0  \\
0 & 0 & -\gamma^2 I 
\end{bmatrix}
\begin{bmatrix} 
 u_1 \\
 y_1 \\
 e_1 
\end{bmatrix}
\right\rangle \\
 \sigma_1(u_1, y_1, e_1) & := \frac{1}{2} \left\langle 
\begin{bmatrix} 
 u_1 \\
 y_1 \\
 e_1 
\end{bmatrix},
\begin{bmatrix}
0 & -I & 0 \\
-I & 0 & I \\
0 & I &0  
\end{bmatrix}
\begin{bmatrix} 
 u_1 \\
 y_1 \\
 e_1 
\end{bmatrix}
\right\rangle.
\end{align*}
Then $G \|_{e_2=0} \Sigma_2$ has $\Ltwo$-gain less than or equal to $\gamma$ for all bounded  passive $\Sigma_2$ if and only if
$\sigma_0(u_1, y_1, e_1) \leq 0$ for all  $(u_1, y_1, e_1) \in \HH$  with $\sigma_1(u_1, y_1, e_1) \geq 0$.
This is equivalent, via the S-procedure lossless theorem, to the existence of $\mu \geq 0$ such that
\[
 \sigma_0(u_1, y_1, e_1) + \mu \sigma_1(u_1, y_1, e_1) \leq 0, \quad \forall (u_1, y_1, e_1) \in \HH.
\]
This implies that
\[
\| y_1 \|_2^2 - \gamma^2  \|e_1\|_2^2  - \mu \langle u_1, y_1 \rangle+  \mu \langle e_1, y_1 \rangle\leq 0, \quad \forall e_1 \in \Ltwo ,
\]
and thus in the subset $\{(u_1, y_1, 0) \in \Ltwo \mid y_1 = G (u_1) \} \subset \HH$, this yields
\[
\mu \langle u_1, y_1 \rangle \geq  \| y_1 \|_2^2, \quad \forall u_1 \in \Ltwo,
\]
i.e., $G$ is output strictly passive. Consequently, $\Sigma_1$ is output strictly passive.

\begin{remark}
Note that by \cite[Prop. 3.1.14]{Sch17} the previous converse passivity theorems extend to the same converse passivity statements for {\it state space systems} that are {\it reachable} from a ground state $x^*$ for which the input-output map $\Delta$ defined by the state space system satisfies the conditions of Theorems \ref{thm: passive}, \ref{thm: passive2}, \ref{Theorem6}.
\end{remark}

Especially the last version of the converse passivity theorem presented in Theorem \ref{Theorem6} is crucial for applications. It implies that closed-loop stability (in the sense of finite $\Ltwo$-gain) of a system
interacting with an unknown, but passive, environment can only be guaranteed if the system seen from the interaction port with the environment is
output strictly passive. This has obvious implications in robotics, where the given system is the controlled robot, interacting with its unknown but physical (and thus typically passive) environment. It is also of importance in the analysis and control of
{\it reduced-order} models, in case the neglected dynamics can be regarded as a passive feedback loop for the reduced-order model. An illustration of this main idea, in a very simple and linear context, is provided in the following example with a robotics flavor.

\begin{example}
Consider an actuated mass
\[
m\dot{v}  =  - dv + u_1, 
\]
where $v$ is the velocity of the mass, $m > 0$ its mass parameter, $d$ the possibly negative `damping' parameter, $u_1$ the external force, and $y_1=v$ the output. Clearly the system is output strictly passive if and only if $d > 0$. Consider an unknown environment modeled by a spring system with spring constant $k>0$ given as
\[
\dot{q}  =  - sq + u_2,  
\]
where $q$ is the extension of the spring, $u_2$ an input velocity and $sq $ a drag velocity (proportional to the spring force $kq$, and thus to $q$). The spring system with output $y_2=kq$ is passive if and only if $s \geq 0$.
The interconnection $u_1=-y_2 +e_1, u_2=y_1$ of the mass system (for arbitrary $d \in \mathbb{R}$) with the spring system results in the closed-loop system
\[
\begin{bmatrix} \dot{v} \\ \dot{q} \end{bmatrix} = 
\begin{bmatrix} -\frac{d}{m} & -\frac{1}{m} \\ 1 & -s \end{bmatrix}
\begin{bmatrix} v \\ q \end{bmatrix} + \begin{bmatrix} \frac{e_1}{m} \\ 0 \end{bmatrix},
\]
with $e_1$ an external force. This system has $\Ltwo$-gain $\leq \gamma$ for some $\gamma >0$ iff the system for $e_1=0$ is
asymptotically stable, which is the case iff $\frac{d}{m} + s > 0$. Hence the closed-loop system has $\Ltwo$-gain $\leq \gamma$ for some
$\gamma >0$ if and only $\frac{d}{m} + s > 0$ for all $s\geq 0$, or equivalently, iff $d >0$, i.e., the mass system is output strictly passive.
%
\end{example}

\section{The converse of the small-gain theorem}

Using similar reasoning as in the passivity case we provide in this section two versions of the {\it converse small-gain theorem}. These results extend the well-known necessity of the small-gain condition for linear systems based on transfer function analysis;
see e.g. \cite{ZDG}. The necessity of the small-gain condition is crucial in robust control theory based on modeling the uncertainty in the `plant'
system by a feedback loop with an unknown system, with magnitude bounded by its $\Ltwo$-gain; see e.g. \cite{ZDG} for the linear case and \cite{Sch17}
(and references therein) for the nonlinear case.

\begin{thm} \label{thm: smallgain} Given $\Sigma_1 = \Delta_2 G \Delta_1$ and $\alpha >0$, where $G$ is linear time-invariant, $\Delta_i$ is linear
  and invertible with $\Ltwo$-gain $= 1$, then for some $\gamma >0$, the closed-loop system $\Sigma_1 \| \Sigma_2$ has $\Ltwo$-gain $\leq \gamma$ for
  all $\Sigma_2$ with $\Ltwo$-gain $\leq \alpha$ if and only if $\Sigma_1$ has $\Ltwo$-gain $ < \frac{1}{\alpha}$.
\end{thm}

\noindent {\bf Proof} \, Sufficiency is well known in the literature. In order to show necessity, note that by the theory of multipliers~\cite[Section
3.5]{GreLim95}, $\Sigma_1 \| \Sigma_2$ having $\Ltwo$-gain $\leq \gamma$ is equivalent to $G \| \Delta_1\Sigma_2\Delta_2$ having $\Ltwo$-gain
$\leq \gamma$, which in turn is equivalent to $G \| \Sigma_2$ having $\Ltwo$-gain $\leq \gamma$, for all  $\Sigma_2$ with $\Ltwo$-gain
$\leq \alpha$. Define
\[
\HH := \{(u_1, u_2, e_1, e_2) \in \Ltwo\ |\ u_2 = e_2 + G (u_1)\}
\]
and the quadratic forms $\sigma_i : \HH
\to \Real$, $i = 0, 1$, as 
\begin{align*}
 \sigma_0(u_1, u_2, e_1, e_2) & := \left\langle 
\begin{bmatrix} 
 u_1 \\
 u_2 \\
 e_1 \\
 e_2
\end{bmatrix},
\begin{bmatrix}
I & 0 & 0 & 0 \\
0 & I & 0 & 0 \\
0 & 0 & -\gamma^2 I & 0 \\
0 & 0 & 0 & -\gamma^2 I
\end{bmatrix}
\begin{bmatrix} 
 u_1 \\
 u_2 \\
 e_1 \\
 e_2
\end{bmatrix}
\right\rangle \\
 \sigma_1(u_1, u_2, e_1, e_2) & := \left\langle 
\begin{bmatrix} 
 u_1 \\
 u_2 \\
 e_1 \\
 e_2
\end{bmatrix},
\begin{bmatrix}
-I & 0 & I & 0 \\
0 & \alpha^2 I & 0 & 0 \\
I & 0 & -I & 0 \\
0 & 0 & 0 & 0
\end{bmatrix}
\begin{bmatrix} 
 u_1 \\
 u_2 \\
 e_1 \\
 e_2
\end{bmatrix}
\right\rangle.
\end{align*}
Then $G \| \Sigma_2$ has $\Ltwo$-gain $\leq \gamma$ for all $\Sigma_2$ with $\Ltwo$-gain $\leq \alpha$ if and only if
\[
 \sigma_0(u_1, u_2, e_1, e_2) \leq 0
 \]
for all  $(u_1, u_2, e_1, e_2) \in \HH$ such that $\sigma_1(u_1, u_2, e_1, e_2) \geq 0$.
This is equivalent, via the S-procedure lossless theorem, to the existence of $\mu \geq 0$ such that
\begin{align*}
 \sigma_0(u_1, u_2, e_1, e_2) + \mu \sigma_1(u_1, u_2, & e_1, e_2) \leq 0, \\
 & \forall (u_1, u_2, e_1, e_2) \in \HH.
\end{align*}
In the subset $\{(u_1, u_2, 0, 0) \in \Ltwo\ |\ u_2 = G (u_1)\} \subset \HH$, this implies that
\[
\|u_1\|_2^2 + \|Gu_1\|_2^2 - \mu \|u_1\|_2^2 + \mu \alpha^2 \|G (u_1)\|_2^2 \leq 0, \; \forall u_1 \in \Ltwo.
\]
It is obvious from the inequality above that $\mu \neq 0$, and hence $\mu >0$. Thus
\[
 \mu \alpha^2 \|G (u_1)\|_2^2 < \mu \|u_1\|_2^2, \quad \forall u_1 \in \Ltwo \mbox{ with } u_1 \neq 0,
\]
and hence $\|G (u_1)\|_2^2 < \frac{1}{\alpha^2} \|u_1\|_2^2$, showing that $G$, and hence $\Sigma_1$, has $\Ltwo$-gain $ < \frac{1}{\alpha}$.

\medskip

In analogy with Theorem \ref{Theorem6} we formulate the following alternative version for the case $e_2=0$.

\begin{thm}
  Given $\Sigma_1 = \Delta_2 G \Delta_1$ and $\alpha >0$, where $G$ is linear time-invariant, $\Delta_i$ is linear and invertible with $\Ltwo$-gain
  $= 1$, then there exists $\gamma$ such that $\Sigma_1 \|_{e_2=0} \Sigma_2$ has $\Ltwo$-gain $\leq \gamma$ from $e_1$ to $y_1$ for all $\Sigma_2$
  with $\Ltwo$-gain $\leq \alpha$ if and only if $\Sigma_1$ has $\Ltwo$-gain $< \frac{1}{\alpha}$.
\end{thm}

\noindent {\bf Proof} \, Sufficiency is clear.  For necessity, note that as in Theorem~\ref{thm: smallgain}, the hypothesis is equivalent to
$G \|_{e_2=0} \Sigma_2$ having $\Ltwo$-gain $\leq \gamma$ for all $\Sigma_2$ with $\Ltwo$-gain $\leq \alpha$. Define
\[
\HH := \{(u_1, y_1, e_1) \in \Ltwo\ |\ y_1 = G (u_1)\}
\]
and the quadratic forms $\sigma_i : \HH
\to \Real$, $i = 0, 1$, as 
\begin{align*}
 \sigma_0(u_1, y_1, e_1) & := \left\langle 
\begin{bmatrix} 
 u_1 \\
 y_1 \\
 e_1 
\end{bmatrix},
\begin{bmatrix}
0 & 0 & 0  \\
0 & I & 0  \\
0 & 0 & -\gamma^2 I 
\end{bmatrix}
\begin{bmatrix} 
 u_1 \\
 y_1 \\
 e_1 
\end{bmatrix}
\right\rangle ,\\
\sigma_1(u_1, y_1, e_1) & := \left\langle 
\begin{bmatrix} 
 u_1 \\
 y_1 \\
 e_1 
\end{bmatrix},
\begin{bmatrix}
-I & 0 & I  \\
0 & \alpha^2 I & 0  \\
I & 0 & -I 
\end{bmatrix}
\begin{bmatrix} 
 u_1 \\
 y_1 \\
 e_1
\end{bmatrix}
\right\rangle.
\end{align*}
Then $\Ltwo$-gain $\leq \gamma$ of $G \|_{e_2=0} \Sigma_2$ for all $\Sigma_2$ with $\Ltwo$-gain $\leq \alpha$ amounts to
\begin{align*}
 \sigma_0(u_1, y_1, e_1) \leq 0 \quad \forall & (u_1, y_1, e_1) \in \HH \\
 & \text{such that} \quad \sigma_1(u_1, y_1, e_1) \geq 0.
\end{align*}
This is equivalent, via the S-procedure lossless theorem, to the existence of $\mu \geq 0$ such that
\[
 \sigma_0(u_1, y_1, e_1) + \mu \sigma_1(u_1, y_1, e_1) \leq 0, \quad \forall (u_1, y_1, e_1) \in \HH.
\]
This implies that
\[
\begin{array}{l}
\| y_1 \|_2^2 - \gamma^2  \|e_1\|_2^2  \\[2mm]
+ \mu \left(u_1(-u_1 + e_1)+  \alpha^2 \|y_1\|_2^2 + e_1(u_1 -e_1) \right) \leq 0
\end{array}
\]
for all $e_1 \in \Ltwo$. Thus in the subset $\{(u_1, y_1, 0) \in \Ltwo\ |\ y_1 = G (u_1) \} \subset \HH$ this yields
\[
\|y_1\|^2 + \mu \alpha^2 \| y_1 \|_2^2 - \mu \|u_1\|^2 \leq 0, \quad  \forall u_1 \in \Ltwo.
\]
This implies $\mu \neq 0$ and thus $\mu >0$, and hence by dividing by $\mu$ it follows that $G$, and hence $\Sigma_1$, has $\Ltwo$-gain $< \frac{1}{\alpha}$. 

\section{Conclusions}

We proved (different versions of the) converse passivity and small-gain theorems for certain linear time-varying systems interconnected in feedback
with nonlinear systems by making crucial use of the S-procedure lossless theorem. Such converse results are fundamental in the control of systems
interacting with unknown environments (e.g., in robotics), and in robust control theory (modeling uncertainty in the to-be-controlled system
by unknown feedback loops). Surprisingly, a full state space version of these results seems to be non-trivial (see \cite{Str15} for partial
results). We also refer to the discussion in \cite{Str15} for further generalizations of the converse passivity theorem; in particular the
quantification of closed-loop stability under interaction with an unknown environment that is allowed to be active in a constrained manner. This is
closely related to the well-known fact that `lack of passivity' of the second system may be `compensated by' excess of passivity, of the first
system; cf. \cite[Theorem 2.2.18]{Sch17} and various work on passivity indices, see e.g. \cite{bao}. Future work also involves seeking converse
results for network-interconnected systems.

\medskip

{\bf Acknowledgement}
The authors are grateful to Carsten Scherer for useful discussions that led to improvements of the paper.

\end{document}